\documentclass[reqno,11pt]{amsart}

\usepackage[all]{xy}
\usepackage{amsthm,array,amssymb,amscd,amsfonts,latexsym}

\DeclareFontFamily{OT1}{wncyr}{\hyphenchar\font45 }
\DeclareFontShape{OT1}{wncyr}{m}{n}{%
   <5> <6> <7> <8> <9> gen * wncyr
   <10> <10.95> <12> <14.4> <17.28> <20.74>  <24.88>wncyr10}{}
\DeclareFontShape{OT1}{wncyr}{m}{it}{%
   <5> <6> <7> <8> <9> gen * wncyi
   <10> <10.95> <12> <14.4> <17.28> <20.74> <24.88> wncyi10}{}
\DeclareFontShape{OT1}{wncyr}{m}{sc}{%
   <5> <6> <7> <8> <9> <10> <10.95> <12> <14.4>
   <17.28> <20.74> <24.88>wncysc10}{}
\DeclareFontShape{OT1}{wncyr}{b}{n}{%
   <5> <6> <7> <8> <9> gen * wncyb
   <10> <10.95> <12> <14.4> <17.28> <20.74> <24.88>wncyb10}{}
\input cyracc.def
\def\rus{\usefont{OT1}{wncyr}{m}{n}\cyracc\fontsize{9}{11pt}\selectfont}

\DeclareMathSizes{9}{9}{7}{5} 

\newtheorem{thm}{Theorem}
\newtheorem{lem}{Lemma}
\newtheorem{cor}{Corollary}
\newtheorem{problem}{Problem}

\theoremstyle{definition}
\newtheorem{defn}{Definition}
\newtheorem{remark}{Remark}
\newtheorem{example}{Example}

\begin{document}

\newcommand{\Span}{\operatorname{Span}}
\newcommand{\Symp}{\mbox{\boldmath$\rm Sp$}}
\newcommand{\g}{\mathfrak{g}}
\newcommand{\el}{\mathfrak{l}}
\newcommand{\lt}{\mathfrak{t}}
\newcommand{\lc}{\mathfrak{c}}
\newcommand{\lu}{\mathfrak{u}}
\newcommand{\lr}{\mathfrak{r}}
\newcommand{\Id}{\operatorname{id}}
\newcommand{\id}{\operatorname{id}}
\newcommand{\pr}{\operatorname{pr}}
\newcommand{\Hom}{\operatorname{Hom}}
\newcommand{\sign}{\operatorname{sign}}

\newcommand{\bbA}{{\mathbb A}}
\newcommand{\bbC}{{\mathbb C}}
\newcommand{\bbZ}{{\mathbb Z}}
\newcommand{\bbP}{{\mathbb P}}
\newcommand{\bbQ}{{\mathbb Q}}
\newcommand{\bbG}{{\mathbb G}}
\newcommand{\bbN}{{\mathbb N}}
\newcommand{\bbF}{{\mathbb F}}

\newcommand{\e}{{\lambda}}

\newcommand{\ve}{{\varepsilon}}
\newcommand{\vp}{{\varpi}}

\newcommand{\kbar}{\overline k}

\newcommand{\mo}{\mathopen<}
\newcommand{\mc}{\mathclose>}

\newcommand{\Ker}{\operatorname{Ker}}
\newcommand{\Aut}{\operatorname{Aut}}
\newcommand{\sdp}{\mathbin{{>}\!{\triangleleft}}} 
\newcommand{\Alt}{\operatorname{A}}   
\newcommand{\GL}{\mbox{\boldmath$\rm GL$}}
\newcommand{\PGL}{\mbox{\boldmath$\rm PGL$}}
\newcommand{\SL}{\mbox{\boldmath$\rm SL$}}
\newcommand{\excG}{\mbox{\boldmath$\rm G$}}
\newcommand{\rank}{\operatorname{rank}}
\newcommand{\aut}{\operatorname{Aut}}
\newcommand{\Char}{\operatorname{\rm char\,}} 
\newcommand{\diag}{\operatorname{\rm diag}}
\newcommand{\Gal}{\operatorname{Gal}}
\newcommand{\galois}{\Gal}
\newcommand{\lra}{\longrightarrow}
\newcommand{\SO}{\mbox{\boldmath$\rm SO$}}
\newcommand{\M}{\operatorname{M}}        
\newcommand{\ord}{\mathop{\rm ord}\nolimits}
\newcommand{\Sym}{{\operatorname{S}}}    
\newcommand{\tr}{\operatorname{\rm tr}}
\newcommand{\trace}{\tr}
\newcommand{\ad}{\operatorname{Ad}}

\newcommand{\Res}{\operatorname{Res}}
\newcommand{\Sha}{\mbox{\rus{\fontsize{11}{11pt}\selectfont{SH}}}}
\newcommand{\G}{\mathcal{G}}
\renewcommand{\H}{\mathcal{H}}
\newcommand{\gen}[1]{\langle{#1}\rangle}
\renewcommand{\O}{\mathcal{O}}
\newcommand{\C}{\mathcal{C}}
\newcommand{\Ind}{\operatorname{Ind}}
\newcommand{\End}{\operatorname{End}}
\newcommand{\Spin}{\mbox{\boldmath$\rm Spin$}}
\newcommand{\T}{\mathbf G}
\newcommand{\GT}{\mbox{\boldmath$\rm T$}}
\newcommand{\Inf}{\operatorname{Inf}}
\newcommand{\Tor}{\operatorname{Tor}}
\newcommand{\m}{\mbox{\boldmath$\mu$}}
\newcommand{\Cr}{\operatorname{Cr}}
\newcommand{\Cay}{\operatorname{Cay}}
\newcommand{\sCay}{\operatorname{sCay}}
\newcommand{\Lie}{\operatorname{Lie}}

\newcommand{\A}{{\sf A}}

\newcommand{\D}{{\sf D}}

\newcommand{\Lbd}{{\sf \Lambda}}


\title[Cayley degree]{On the Cayley degree of an algebraic group}
\author{Nicole Lemire}
\address{Department of Mathematics, University of Western Ontario,
London, Ontario N6A 5B7, Canada} \email{nlemire@uwo.ca}

\author{ Vladimir L. Popov}
\address{Steklov Mathematical Institute,
Russian Academy of Sciences, Gubkina 8, Moscow 119991,
Russia} \email{popovvl@orc.ru}

\author{Zinovy Reichstein}
\address{Department of Mathematics, University of British Columbia,
       Vancouver, BC V6T 1Z2, Canada}
\email{reichste@math.ubc.ca}
\urladdr{www.math.ubc.ca/$\stackrel{\sim}{\phantom{.}}$reichst}
\thanks{N. Lemire and Z. Reichstein
were supported in
 part by  NSERC research grants. \\
${\ }$\quad V. L. Popov
 was supported  in part by ETH, Z\"urich
 (Switzerland), Russian grants {\rus RFFI
05--01--00455}, {\rus N{SH}--9969.2006.1}, and a
(granting) program of the Mathematics Branch of the
Russian Academy
of Sciences.}


\begin{abstract} A connected linear algebraic group $G$ is
called a {\em Cayley group} if
the Lie algebra of $G$ endowed with the adjoint
$G$-action and the
group variety of $G$ endowed with the conjugation
$G$-action are birationally $G$-isomorphic. In particular, the
classical Cayley map
\[ X \mapsto (I_n-X)(I_n+X)^{-1} \]
between the special orthogonal group $\SO_n$ and its Lie
algebra ${\mathfrak so}_n$, shows that $\SO_n$ is a
Cayley group. In an earlier paper we
classified the simple Cayley groups defined over an
algebraically closed field of characteristic zero.
Here we consider a new numerical invariant of $G$, {\em
the Cayley degree}, which ``measures" how far $G$ is
from being Cayley, and prove upper bounds on Cayley
degrees of some groups.
\end{abstract}

\maketitle

\section {Introduction}
\label{sect.intro}

Let $G$ be a connected linear algebraic group and let
$\g$ be its Lie algebra.  We say that $G$
is a {\em Cayley group} if there is a birational
isomorphism
\begin{equation} \label{e.Cayley-map}
\varphi \colon G \dasharrow \g
\end{equation}
which is equivariant with respect to the conjugation
action of $G$ on itself and the adjoint action of $G$ on
$\g$; see \cite[Definition 1.5]{lpr}. In particular, the
classical Cayley map~\cite{cayley}
\begin{equation} \label{e.classical}
X \mapsto (I_n-X)(I_n+X)^{-1}
\end{equation}
between the special orthogonal group $\SO_n$ and its Lie
algebra ${\mathfrak so}_n$
shows that $\SO_n$ is a Cayley group. (The same formula
shows that $\Symp_{2n}$ is Cayley as well.) In the
sequel we will always assume that the base field $k$ is
algebraically closed and of characteristic zero.
(Problem~\ref{prob.cay-deg} below is of interest for
arbitrary $k$ but the partial answers we would like to
discuss here require this assumption.)

In 1975 D. Luna~\cite{lu1}, \cite{lu3} asked the
second-named author a question that, in the above
terminology, can be restated as follows: For what $n$ is the
group $\SL_n$ Cayley? In~\cite{lpr} we showed that
$\SL_n$ is Cayley if and only if $n \leqslant 3$ and, more
generally, proved the following classification theorem.

\begin{thm} \label{thm.lpr1}
{\rm(}\cite[Theorem 3.31(a)]{lpr}{\rm)} A connected
simple algebraic group $G$ is Cayley if and only $G$ is
isomorphic to one of the following groups:
 $\SL_2$, $\SL_3$, $\SO_{n}$ $(n \ne 2, 4)$,
$\Symp_{2n}$,  $\PGL_n$ $(n \geqslant 1)$.
\end{thm}

Note that $\SO_n$ is a Cayley group for every $n
\geqslant 1$; we have excluded $\SO_2$ and $\SO_4$ from
the above list because these groups are not simple.

A {\em generalized Cayley map} of $G$ is a rational
$G$-equivariant map $\varphi \colon G \dasharrow \g$, as
in \eqref{e.Cayley-map}, except that instead of
requiring it to be a birational isomorphism, we only
require it to be dominant, see~\cite[Definition
10.9]{lpr}. Every generalized Cayley map of $G$ has
finite degree,
\[ \deg \varphi = [k(G) : \varphi^*(k(\g))]<\infty, \]
where, as usual, $k(X)$ and
$k[X]$ denote respectively the field of rational and the
algebra of regular functions on an irreducible algebraic
variety $X$). A generalized Cayley
map~\eqref{e.Cayley-map} exists for every linear
algebraic group $G$; see \cite[Proposition 10.5]{lpr}.
Hence the following natural number is well defined.

\begin{defn} \label{def.Cayley-degree}
The {\em Cayley degree} $\Cay(G)$ of $G$ is 
the minimal value of $\deg \varphi$, as $\varphi$ ranges
over all generalized Cayley maps of $G$.
\end{defn}

Note that, by definition, $G$ is a Cayley group if and only
if $\Cay(G)=1$. Therefore
Theo\-rem~\ref{thm.lpr1} may be viewed as a first step
toward a solution of the following more general problem.

\begin{problem} \label{prob.cay-deg} {\it Find
the Cayley degrees of connected simple algebraic
groups.}
\end{problem}

We do not have any general methods for proving lower
bounds on the Cayley degree, beyond those provided by
Theorem~\ref{thm.lpr1}; in particular, we do not have an
example of a linear algebraic group $G$ with $\Cay(G)>2$.
Thus in this note we will primarily concentrate on upper
bounds.  Our main results are Theorems~\ref{thm.main1}
and~\ref{thm.main2} below.

\begin{thm} \label{thm.main1}
If $n \geqslant 3$, then $\Cay(\SL_n) \leqslant n-2$.
\end{thm}

Our proof of Theorem~\ref{thm.main1} 
is self-contained. For $n = 3$ this argument gives a
new proof of the fact that $\Cay(\SL_3) = 1$ (i.e.,
$\SL_3$ is a Cayley group), which is simpler than either
of the two proofs in~\cite{lpr}.  For $n = 4$,
Theorem~\ref{thm.main1} implies that $\Cay(\SL_4) = 2$;
see Example~\ref{ex.isogeny4}.

To motivate our second main result, we note that the exceptional
group $\excG_2$ plays a special role in this theory.
While $\excG_2$ is not a Cayley group, it is close to
being one, in the sense that $\excG_2 \times {\bf
G}_m^2$ is Cayley; see~\cite[Theorem 1.31(b)]{lpr}. In
fact, $\excG_2$ is the unique simple group $G$ which is
stably Cayley but is not Cayley; see \cite[Theorems 1.29
and 1.31]{lpr}. (Recall that $G$ is called {\em stably
Cayley} if $G \times {\bf G}_m^r$ is Cayley for some $r
\geqslant 1$.) Theorem~\ref{thm.main2} below shows that
$\excG_2$ is also close to being Cayley in the sense of
having a small Cayley degree.

\begin{thm} \label{thm.main2}
$\Cay(\excG_2) = 2$.
\end{thm}

The rest of this note is structured as follows. In
Section~\ref{sect2} we determine the Cayley degrees
of Spin groups and some groups of type ${\sf A}$. In
Section~\ref{sect3} we prove a lemma that reduces
the computation of the Cayley degree of a reductive
group $G$ to a question about finite group actions. This
lemma is then used as a starting point for the proofs of
Theorems~\ref{thm.main1} and \ref{thm.main2} in
Sections~\ref{sect4} and~\ref{sect5} respectively.
In Section~\ref{sect6} we give a representation
theoretic interpretation of the Cayley degree.

\section{First examples}
\label{sect2}

\begin{lem} \label{lem.isogeny} 

\begin{enumerate}
\item[\rm(a)]
Let $\pi \colon G\to H$ be an isogeny
between connected linear algebraic groups and
let $d$ be the order of its kernel.

\begin{enumerate}
\item[\rm(a$_1$)] Then
\[ \Cay(G)\leqslant d \cdot \Cay(H) \, . \]

\smallskip
\item[\rm(a$_2$)] If $G$ is not Cayley but $H$ is
Cayley, and $d = 2$, then $\Cay(G) =2$.
\end{enumerate}

\smallskip
\item[\rm(b)] Let $\varphi_i$ be a generalized Cayley
map of a connected linear algebraic group $G_i$, where
$i=1,\ldots, n$. Then
$\varphi_1\times\ldots\times\varphi_n$ is a generalized
Cayley map of $G_1\times\ldots\times G_n$, and
$$\deg \,(\varphi_1\times\ldots\times\varphi_n)=
\deg \,\varphi_1\ldots \deg \,\varphi_n.$$
\end{enumerate}
\end{lem}

\begin{proof}
(a$_1$) The groups $G$ and $H$ have the same Lie algebra
$\g$. Let $\varphi\colon H\dasharrow \g$ be a
generalized Cayley map of $H$. Since $\Ker \pi$ is a
finite central subgroup of $G$ and $\deg \pi=d$, the
composition $\varphi \circ\pi\colon G\dasharrow \g$ is a
generalized Cayley map of $G$. Its degree is $d \cdot
\deg \,\varphi$, and part (a$_1$) follows.

(a$_2$) Since $G$ is not Cayley, we have $\Cay(G)
\geqslant 2$. The opposite inequality follows from part
(a$_1$).

Part (b) follows from the interpretation of degree
of a rational map as the number of points in a general fiber.
\end{proof}

From (b) and Definition~\ref{def.Cayley-degree} we
obtain the following upper bound.

\begin{cor} \label{cor.ineq}
$\Cay(G_1\times\ldots\times G_n)\leqslant
\Cay(G_1) \cdot \ldots \Cay (G_n)$.
\end{cor}
The following example shows that, in general, equality does not hold.

\begin{example} Since
$\Cay(\excG_2) \ge 2$ by Theorem~\ref{thm.lpr1}, but
$\Cay({\bf G}_2\times{\bf G}_m^2)=1$ (see \cite[Theorem 1.31]{lpr}),
we see that
\[ \Cay({\bf G}_2\times{\bf G}_m^2) <
\Cay(\excG_2) \cdot \Cay({\bf G}_m^2) \, . \] (In fact,
the right hand side of this inequality is equal to $2$,
because $\Cay(\excG_2) = 2$ by Theorem~\ref{thm.main1}
and $\Cay(\bbG_m^2) = 1$; see~\cite[Example 1.21]{lpr}.)
\end{example}

\begin{example} \label{ex.isogeny2} (see \cite[p. 962]{lpr})
The groups
\[ \Spin_2 \simeq \bbG_m \, ,\quad
   \Spin_3 \simeq \SL_2 \, , \quad
   \Spin_4 \simeq \SL_2 \times \SL_2 \, , \quad
   \Spin_5 \simeq \Symp_4 \,  \]
are easily seen to be Cayley. On the other hand,
$\Spin_n$ is not Cayley if $n \geqslant 6$. Since
$\SO_n$ is Cayley for every $n$, applying
Lemma~\ref{lem.isogeny}(b) to the natural $2$-sheeted
isogeny $\Spin_n \to \SO_n$ (where $n \geqslant 6$), we
obtain
\begin{equation}\label{cspin}
\Cay ({\bf Spin}_n)=
\begin{cases} 2&\text{for $n\geqslant 6$},\\
1&\text{for $n\leqslant 5$}.
\end{cases}
\end{equation}
\end{example}

\begin{example} \label{ex.isogeny3}
Since $\PGL_n$ is a Cayley group for every $n \geqslant
1$, Lemma~\ref{lem.isogeny}, applied to the natural
isogeny $\SL_n/\m_d =: G \to H:=\PGL_n$ yields
\begin{equation} \label{e.typeA}
\Cay(\SL_n/\m_d)\leqslant n/d \, .
\end{equation}
In particular,
\[ \Cay(\SL_{2d}/\m_d)=\begin{cases}
2&\text{for $d\geqslant 3$},\\
\label{sl2m} 1&\text{for $d\leqslant 2$.}
\end{cases} \]
Note also that setting $d = 1$ in \eqref{e.typeA}
yields $\Cay(\SL_n)\leqslant n$. Theorem~\ref{thm.main1}
strengthens this bound.
\end{example}

\section{The maximal torus}
\label{sect3}

In this section we reduce the problem of finding
$\Cay(G)$ for a connected reductive group $G$, to a
question about finite group actions.

\begin{lem} \label{lem.torus}
Let $G$ be a connected linear algebraic group, let $T$
be its maximal torus, let $C$ and $N$ be the centralizer
and normalizer of $T$ in $G$ respectively, and let $W :=
N/C$ be the Weyl group. Denote the Lie algebras of $G$,
$T$, and $C$ by $\g$, $\lt$, and $\lc$, respectively.
\begin{enumerate}
\item[\rm(a)] Then
\begin{equation}\label{eqv}
\Cay(G)= \min_{\psi}\,\deg \psi,
\end{equation}
where $\psi$ ranges over all dominant rational
$N$-equivariant maps $C\dasharrow \lc$.
\item[\rm(b)] Moreover, if $G$ is reductive, then
\eqref{eqv} holds, where $\psi$ ranges over all
$W$-equivariant dominant rational maps $T \dasharrow \lt$.
\end{enumerate}
\end{lem}

\begin{proof} Recall that  $ G \simeq G\times^N C$ and
$\g \simeq G\times^N \lc$, where $\simeq$ stands for a
birational isomorphism of $G$-varieties.
Moreover, if $\varphi \colon G\times^N C \dasharrow G
\times ^N\lc$ is a dominant rational $G$-map, then
$\psi:=\varphi|_C: C\dasharrow\lc$ is a dominant
rational
$N$-map and $\varphi^{-1}(x)=
\psi^{-1}(x)$ for a general point $x\in \lc$;
see~\cite[Lemma 2.17]{lpr}. Hence
\begin{equation}\label{deg-restr}
\deg\varphi =|{\varphi^{-1}(x)}|=
|\psi^{-1}(x)| = \deg \psi.
\end{equation}
Thus we have a degree preserving bijection between
generalized Cayley maps of $G$ and dominant rational
$N$-equivariant maps $C\dasharrow \lc$. This immediately
implies (a).
If $G$ is reduc\-tive, then $C = T$, $\lc = \lt$, and
the $N$-actions on $C$ and $\lc$ descend to the
$W$-actions (since $T$, being commutative, acts
trivially). Hence part (b) follows from part~(a).
\end{proof}

\begin{cor}\label{degree:} Let $\varphi$ be a
generalized Cayley map of a connected reductive group
$G$. Then
$\deg \,\varphi=[k(G)^G:\varphi^*(k(\g)^G)]$.
\end{cor}

\begin{proof} We will continue to use the notations of
Lemma~\ref{lem.torus} and set $\psi:=\varphi|_T$. Since $W$
is a finite group acting on $T$ and $\lt$ faithfully, we have
$[k(T):k(T)^W]=\!|W|$ and $[k(\lt):k(\lt)^W]\!=\!|W|$.
From this we deduce that $\deg \, \psi:=\![k(T):\psi^*(k(\lt))]\!=
\![k(T)^W:\psi^*(k(\lt)^W)]$. Since we have
$[k(T)^W:\psi^*(k(\lt)^W)]=[k(G)^G:\varphi^*(k(\g)^G)]$,
see \cite[Theorem (1.7.5)]{popov}, \cite[(3.4)]{lpr},
the claim now follows from \eqref{deg-restr}.
\end{proof}

\begin{remark}
If $\varphi$ is a morphism, Corollary~\ref{degree:}
can be deduced from \cite[Lemme Fondamental]{lu3}. For
certain particular morphisms $\varphi$,  a proof can be
found in \cite[Corollary (3.3)]{komi}.
\end{remark}

\section{Proof of Theorem~\ref{thm.main1}}
\label{sect4}

By Lemma~\ref{lem.torus} it suffices to construct a
dominant rational $W = \Sym_n$-equivariant map between
the maximal torus $T$ in $\SL_n$ and its Lie algebra
$\lt$.

To keep the notation clear in the construction to
follow, we will work with two copies of the affine space
$\bbA^n$, with the same natural (permutation) action of
$\Sym_n$. We will denote one by $\bbA_x^n$ and the other
by $\bbA_y^n$ and use the variables $x_1, \dots, x_n$ and,
respectively, $y_1, \dots, y_n$ as standard coordinate
functions on $\bbA_x^n$ and $\bbA_y^n$.
We will now embed $\lt$ and, respectively, $T$
into $\bbA_x^n$ and $\bbA_y^n$ as the following
$\Sym_n$-invariant subvarieties:
\begin{gather*}
\lt = \{ (a_1, \dots, a_n)\in \bbA_x^n \mid a_1
+ \dots + a_n
= 0 \}, \\
T = \{ (b_1, \dots, b_n)\in \bbA_y^n \mid b_1 \dots b_n = 1 \}.
\end{gather*}
Consider the mutually inverse $\Sym_n$-equivariant
rational maps $\varphi \colon \bbA_x^n \to \bbA_y^n$
and $\psi
\colon \bbA_y^n \to \bbA_x^n$ given by
\[
\varphi
:= \biggl(\frac{x_1 + 1}{x_1}, \ldots, \frac{x_n +
1}{x_n}\biggr) \quad \text{and} \quad
\psi:= \biggl(\frac{1}{y_1-1}, \ldots, \frac{1}{y_n
-1}\biggr).
\]
These maps give rise to a (biregular) isomorphism
between the open subsets
\[ U_x := \{ (a_1, \dots, a_n)\in\bbA_x^n \mid a_1 \ldots
a_n \ne 0 \}
\] and
\[ U_y := \{ (b_1, \dots, b_n)\in \bbA_y^n \mid (b_1 - 1)
\dots (b_n -1) \ne 0 \}
\]
in $\bbA_x^n$ and $\bbA_y^n$ respectively.  Substituting
$y_i = \dfrac{x_i +1}{x_i}$ into the equation $y_1 \dots
y_n -1 = 0$ of $T$, we see that $\psi(T \cap U_y) = X
\cap U_x$, where $X$ is the hypersurface in $\bbA_x^n$
cut out by the equation
\[ f(x_1, \dots, x_n) := (x_1 + 1) \dots (x_n + 1)
- x_1 \dots x_n = 0. \] Since $X \cap U_x$ is isomorphic
to $T \cap U_y$ (which is irreducible) and $X$ does not
contain any of the $n$ components $\{ x_i = 0 \}$ of the
complement of $U_x$, we conclude that $X$ is irreducible
$\Sym_n$-invariant hypersurface in $\bbA_x^n$. Hence $f$
is a power of an irreducible polynomial. Since
$\deg \, f(1,\ldots,1, x_i, 1,\ldots , 1)=1$ for every $i$,
we conclude that in fact $f$ is irreducible. As
$\deg \, f=n-1$, this implies that $X$ is a hypersurface of
degree $n-1$. By our construction $X$ is birationally
isomorphic to $T$ (via $\varphi$), as an
$\Sym_n$-variety.

Let $\pi$ be the projection $X \dasharrow \lt$ from a
point $\mbox{\boldmath$a$} = (a, \dots, a) \in
\bbA_x^n$. That is, for any point $\mbox{\boldmath$b$}
\in X$, $\mbox{\boldmath$b$} \neq \mbox{\boldmath$a$}$,
the point $\pi(\mbox{\boldmath$b$})$ is the intersection
point of the line passing through $\mbox{\boldmath$a$}$
and $\mbox{\boldmath$b$}$ with the hyperplane $\lt
\subset \bbA_x^n$. Moreover, we choose
$\mbox{\boldmath$a$}$ so that it lies on $X$. Note that
this automatically means that it does not lie in $\lt$.
Indeed, since zero does not satisfy
the equation
\[ f(a, \dots, a) = (1+ a)^n - a^n = 0, \]
if $\mbox{\boldmath$a$}\in X$, then
$\mbox{\boldmath$a$}$ cannot lie in $\lt$. Since our
base field $k$ is algebraically closed and of
characteristic zero, such an $a$ exists for every $n
\geqslant 2$. Note that $\pi$ is well-defined, unless
$X$ is a hyperplane parallel to $\lt$.
Since $\deg X = n-1$, it is not a hyperplane for every
$n \geqslant 3$. Thus $\pi$ is well-defined for every $n
\geqslant 3$. Note also that since $\mbox{\boldmath$a$}$
is fixed by $\Sym_n$, the map $\pi$ is
$\Sym_n$-equivariant.

We claim that $\pi \colon X \dasharrow \lt$ is dominant.
Since $\pi$ is a projection map from a point on a
hypersurface $X$, and $\deg X = n-1$, this claim implies
that $\deg\,\pi = n-2$. Composing $\pi$ with a
birational isomorphism $\psi \colon T \dasharrow X$, we
obtain an $\Sym_n$-equivariant dominant rational map $T
\dasharrow \lt$ of degree $n-2$, and
Theorem~\ref{thm.main1} is proved.

It remains to show that $\pi$ is dominant. Assume the
contrary. Let $X_0$ be the closure of the image of $\pi$
in $\lt$. Then $X$ is the cone over $X_0$ centered at
${\bf a}$. Since, as we remarked above, $X$ is not a
hyperplane (we are assuming throughout that $n \geqslant
3$), $X$ has to be singular at $\mbox{\boldmath$a$}$.
Consequently, $a$ satisfies the system of equations
\[ \left \{ \begin{array}{l}
 f(\mbox{\boldmath$a$}) = (1+ a)^n - a^n = 0 \, , \\
\dfrac{\partial f}{\partial x_1}(\mbox{\boldmath$a$}) =
 (1+ a)^{n-1} - a^{n-1} = 0 \, . \end{array} \right.\]
But this system has no solutions, a contradiction.
Theorem~\ref{thm.main1} is now
proved. \qed

\begin{example} \label{ex.isogeny4}
By Theorem~\ref{thm.main1}, $\Cay(\SL_4) \leqslant 2$.
Equivalently, $\Cay(\SL_4) = 2$; indeed, we know that
$\Cay(\SL_4) \ne 1$, i.e., $\SL_4$ is not a Cayley group
by Theorem~\ref{thm.lpr1}.

Since $\SL_4 /\m_2 \simeq \SO_4$ is Cayley, the equality
$\Cay(\SL_4) = 2$ can also be obtained by applying
Lemma~\ref{lem.isogeny}(b) to the isogeny $\SL_4 \to
\SL_4/\m_2$.  Alternatively, since $\SL_4\simeq\Spin_6$,
the equality $\Cay(\SL_4) = 2$ is a special case
of~\eqref{cspin}.
\end{example}

\section{Proof of Theorem~\ref{thm.main2}}
\label{sect5}

First recall that $\excG_2$ is not Cayley
(see~Theorem~\ref{thm.lpr1}) and hence $\Cay(\excG_2)
\geqslant 2$. Thus we only need to prove the opposite
inequality. By Lemma~\ref{lem.torus} it suffices to
construct a $W$-equivariant dominant rational map $T
\dasharrow \lt$ of degree $2$, where $T$ is a maximal
torus of $\excG_2$, $\lt$ is the Lie algebra of $T$, and
$W$ is the Weyl group.

Recall that 
$W$ is isomorphic to $\Sym_3 \times \bbZ/2 \bbZ$. Once
again, we consider two copies of the $3$-dimensional
affine space, $\bbA_x^3$ and $\bbA_y^3$, with the
following $W$-actions. The symmetric group $\Sym_3$ acts
on both copies in the natural way (by permuting the
coordinates).  The nontrivial element of $\bbZ/2 \bbZ$
acts on $\bbA_x^3$ by
\[ (a_1, a_2, a_3) \mapsto (- a_1, - a_2, - a_3), \]
and on $\bbA_y^3$ by
\[ (b_1, b_2, b_3) \mapsto \biggl(\frac{1}{b_1},
\frac{1}{b_2}, \frac{1}{b_3}\biggr). \]
 We may (and shall) embed $\lt$ and $T$
 into $\bbA_x^3$ and $\bbA_y^3$, respectively,
 as the following $W$-invariant subvarieties:
 \begin{gather*}
 \lt = \{ (a_1, a_2 , a_3)\in \bbA_x^3 \mid a_1
 + a_2 + a_3
 = 0 \}, \\
T = \{ (b_1, b_2, b_3)\in \bbA_y^3 \mid b_1 b_2 b_3
= 1 \}.
\end{gather*}

We now consider the mutually inverse $W$-equivariant
rational maps $\varphi \colon \bbA_x^3 \to \bbA_y^3$ and
$\psi \colon \bbA_y^3 \to \bbA_x^3$ given by
\[
\varphi:= \biggl(\frac{x_1 -1}{x_1 + 1}, \frac{x_2 -
1}{x_2 + 1}, \frac{x_3 - 1}{x_3 + 1}\biggr) \quad
\text{and} \quad
\psi:= \biggl(-\frac{y_1 + 1}{y_1 - 1}, -\frac{y_2 +
1}{y_2 - 1}, -\frac{y_3 + 1}{y_3 - 1}\biggr).
\]
These maps give rise to a $W$-equivariant isomorphism
between the open subsets
\[ U_x := \{ (a_1, a_2, a_3)\in  \bbA_x^3 \mid
(a_1 + 1) (a_2 + 1) (a_3 + 1) \ne 0 \} 
\] and
\[ U_y := \{ (b_1, b_2, b_3)\in  \bbA_y^3 \mid (b_1 - 1) (b_2 - 1)
(b_3 - 1) \ne 0 \}
\]
in  $\bbA_x^3$ and  $\bbA_y^3$, respectively.
Substituting $y_i = \dfrac{x_i -1}{x_i + 1}$ into the
equation $y_1 y_2 y_3 = 1$ of $T$, we see that $\psi(T
\cap U_y) = X \cap U_x$, where $X$ is the $W$-invariant
quadric surface in $\bbA_x^3$ 
defined by the equation
\[ x_1x_2 + x_2x_3 + x_1x_3 + 1 = 0. \]
Composing the $W$-equivariant birational isomorphism
$\psi \colon T \dasharrow \bbA_x^3$ with the
$W$-invariant linear projection $\alpha \colon X \to
\lt$ given by
\[ \alpha
:= \biggl(x_1 - \frac{x_1 + x_2 + x_3}{3}, x_2 -
\frac{x_1 + x_2 + x_3}{3}, x_3 - \frac{x_1 + x_2 +
x_3}{3}\biggr),
\] we obtain a desired $W$-equivariant rational map 
$\alpha \circ\psi \colon T \dasharrow \lt$ of degree $2$. \qed

\begin{remark} The proofs of Theorems~\ref{thm.main1} and~\ref{thm.main2}
proceed along similar lines: we begin by defining a
birational isomorphism $\psi$ between $T$ and a
hypersurface $X$, then project $X$ onto $\lt$. Note,
however, that the projections $\pi$ (in the proof of
Theorem~\ref{thm.main1}) and $\alpha$ (in the proof of
Theorem~\ref{thm.main2}) are different in the following
sense: $\pi$ is a projection from a point on $X$, and
$\alpha$ is a linear projection ($\alpha$ may also be
viewed as a  projection from a point at infinity, which
does not lie on $X$). Note that $\alpha$ cannot be
replaced by a projection from a point of $X$, since $X$
has no $W$-equivariant points (and also because
otherwise $\alpha$ would have degree $1$ and our
argument would show that $\excG_2$ is a Cayley group,
which we know to be false).
\end{remark}

\begin{remark}
The formula for $\varphi$ is somewhat similar to the
formula for the classical Cayley
map~\eqref{e.classical}. Note, however, that we cannot
replace $\dfrac{x_1 -1}{x_1 + 1}$, $\dfrac{x_2 -1}{x_2 +
1}$, etc. by $\dfrac{1 - x_1}{x_1 + 1}$, $\dfrac{1 -
x_2}{x_2 + 1}$, etc. in the definition of $\varphi$. If
we do this, then, setting $\psi = \varphi^{-1}$, we see
that the image of $T$ under $\psi$ becomes the cubic
$x_1 x_2 x_3 + x_1 + x_2 + x_3 = 0$, rather than the
quadric $x_1 x_2 + x_2 x_3 + x_1 x_3 + 1 = 0$, and the
above argument gives a generalized Cayley map of degree
$3$, rather than $2$.
\end{remark}

\section{A representation theoretic approach}
\label{sect6}

In conclusion we outline a representation theoretic
approach to determining the Cayley degree of an algebraic group.

Let $X$ be an irreducible algebraic variety endowed with
an action of an algebraic group $H$, and let $V$ be a
vector space over $k$ of dimension $\dim X$ endowed with
a linear action of $H$. Then rational dominant $H$-maps
$X\dasharrow V$ are described as follows. Let $M$ be a
submodule of the $H$-module $k(X)$ such that
\begin{enumerate}
\item[\rm (i)] $M$ is isomorphic to the $H$-module
$V^*$, \item[\rm (ii)] $k(X)$ is algebraic over the
subfield $k(M)$ generated by $M$ over $k$.
\end{enumerate}
By (ii), $k(M)/k$ is a purely transcendental extension
of degree $\dim X$. Since $k(V)$ is gene\-ra\-ted over
$k$ by $V^*$, any isomorphism of $H$-modules $V^*\to M$
can be uniquely extended up to an $H$-equivariant
embedding $\iota\!: k(V)\hookrightarrow k(X)$ whose
image is $k(M)$. This embedding determines a rational
dominant $H$-map $\psi\!: X\dasharrow V$ such that
$\psi^*=\iota$. We have \begin{equation}\label{repr}
\deg \,\psi=[k(X):k(M)]. \end{equation} Any dominant
rational $H$-map $X\dasharrow V$ is obtained in this
way.

Now suppose $G$ is a connected reductive linear algebraic group,
$X = T$ is a maximal torus, $V = \lt$ is the Lie algebra of $T$ and
$H = W = N_G(T)/T$ is the Weyl group.  In view of
Lemma~\ref{lem.torus}(b) the above approach relates
generalized Cayley maps of $G$ to the $W$-module structure of
$k(T)$. This connection may be used to prove upper bounds on $\Cay(G)$.

\begin{example} Let $G={\bf G}_2$.
Use the notation of Section~\ref{sect5}. Let $t_i$ be
the restriction of $y_i$ to $T$.
Then $t_1t_2t_3=1$ and $k(T)=k(t_1, t_2)$. Put
\begin{equation}\label{yi}
z_i:=t_i-t_i^{-1}.
\end{equation}

 From the description of the $W$-actions on $T$ and $\lt$
given in Section \ref{sect5} it follows that
\begin{equation}\label{M}
M:=\{\alpha_1z_1+\alpha_2 z_2+\alpha_3 z_3\mid
\alpha_1+\alpha_2+\alpha_3=0,\hskip 2mm \alpha_i\in k\}
\end{equation}
is a submodule of the $W$-module $k(T)$ that is
isomorphic to the $W$-module $\lt^*$. Let
\begin{equation}\label{zi}
s_1:=z_1-z_2, \ s_2:=z_1-z_3
\end{equation}
Then $s_1, s_2$ is a basis of $M$, so $k(M)=k(s_1,
s_2)$. We have $k(t_1, s_1, s_2)=k(T)$ because
$t_2=(t_1^2-1)(t_1^2s_1+
t_1s_2-t_1^3-t_1^2+t_1+1)^{-1}$. It follows from
\eqref{yi}, \eqref{zi} that
\begin{equation}\label{system}
\begin{cases}
-t_2+t_2^{-1}=s_1-t_1+t_1^{-1},\\
t_1t_2-t_1^{-1}t_2^{-1}=s_2-t_1+t_1^{-1}.
\end{cases}
\end{equation}

\noindent Eliminating $t_2$ and $t_2^{-1}$ from
\eqref{system}, we obtain the following equation:
\begin{multline*}
t_1^6-(s_1+s_2)t_1^5+(s_1s_2-2s_1-2s_2-1)t_1^4+
(s_1^2+s_2^2-5)t_1^3\\
+(s_1s_2+2s_1+2s_2+1)t_1^2+ (s_1+s_2+1)t_1+1=0.
\end{multline*}

Thus for the conjugating and adjoint actions of $H:=W$
respectively on $X:=T$ and $V:=\lt$, and for $M$ defined
by \eqref{M}, the above conditions (i), (ii) hold and
$[k(T):k(M)]\leqslant 6$. Hence by \eqref{repr},
\eqref{deg-restr}, and Lemma \ref{lem.torus}, there
exists a generalized Cayley map of $G$ of degree $[k(T):k(M)]$. In
particular, this implies that
$\Cay(\bbG_2)\leqslant\!~6$ (of course, by
Theorem~\ref{thm.main2}, we know that in fact
$\Cay(\bbG_2)=2$). \quad $\square$
\end{example}

\end{document}